\documentclass[a4paper,11pt,reqno]{amsart}
\usepackage{enumerate}
\usepackage{graphicx}
\usepackage{amsfonts}
\usepackage{graphicx}
\usepackage{amsmath}
\usepackage{amssymb}
\newcommand{\RR}{{\mathbb R}}
\newcommand{\CC}{{\mathbb C}}
\newcommand{\N}{{\mathbb N}}
\newcommand{\KK}{{\mathbb K}}

\providecommand{\C}[1]{\mathcal{#1}}
\newcommand{\ee}{\mathcal{E}}
\newcommand{\Dee}{\mathcal{D}}
\DeclareMathOperator{\supp}{supp}
\DeclareMathOperator{\capa}{cap}
\newtheorem{theorem}{Theorem}[section]

\newtheorem{assumption}[theorem]{Assumption}

\newtheorem{coro}[theorem]{Corollary}

\newtheorem{definition}[theorem]{Definition}

\newtheorem{lemma}[theorem]{Lemma}

\newtheorem{prop}[theorem]{Proposition}
\newtheorem{remark}[theorem]{Remark}

\newcommand{\Hmm}[1]{\leavevmode{\marginpar{\tiny%
$\hbox to 0mm{\hspace*{-0.5mm}$\leftarrow$\hss}%
\vcenter{\vrule depth 0.1mm height 0.1mm width \the\marginparwidth}%
\hbox to 0mm{\hss$\rightarrow$\hspace*{-0.5mm}}$\\\relax\raggedright #1}}}
\begin{document}
% -------------------------------------------------------------%
\title[Sch'nol's Theorem]{Sch'nol's Theorem For Strongly Local Forms}
% -------------------------------------------------------------%
\author[A.~Boutet de Monvel]{Anne Boutet de Monvel$^1$}
\author[D.~Lenz]{Daniel Lenz$^2$} \author[P.~Stollmann]{Peter
  Stollmann$^2$} \address{$^1$ IMJ, case 7012, Universit\'e Paris 7, 2
  place Jussieu, 75251 Paris, France} \address{$^2$ Fakult\"at f\"ur
  Mathematik, Technische Universit\"at, 09107 Chemnitz, Germany}
% -------------------------------------------------------------%
% -------------------------------------------------------------%
\begin{abstract}
  We prove a variant of Sch'nol's theorem in a general setting: for
  generators of strongly local Dirichlet forms perturbed by measures.

  As an application, we discuss quantum graphs with $\delta$- or
  Kirchhoff boundary conditions.
\end{abstract}
% -------------------------------------------------------------%
\date{\today} %
\maketitle
% -------------------------------------------------------------%
\begin{center}
\emph{Dedicated to Shmuel Agmon on the occasion of his 85th birthday}
\end{center}

\section*{Introduction}
The behavior of solutions to elliptic partial differential equations and its interplay with spectral properties of the associated partial differential operators is a topic of fundamental interest. Our understanding today is in many aspects based on groundbreaking work by Shmuel Agmon (cf \cite{Agmon-65,Agmon-70,AgmonH-76,Agmon-82,Agmon-84,Agmon-85}) to whom this article is dedicated with great admiration and gratitude.
Here we explore the  well known classical fact that the spectral values of
Schr\"odinger operators $H$ can be characterized in terms of the
existence of appropriate ``generalized eigenfunctions'' or
``eigensolutions''.  One part of this characterization is
\emph{Sch'nol's theorem} stating that existence of an eigensolution of
$Hu=\lambda u$ ``with enough decay'' guarantees $\lambda\in\sigma(H)$.
We refer to the original result \cite{Schnol-57} by Sch'nol from 1957
that was rediscovered by Simon, \cite{Simon-81c}, as well as the
discussion in \cite{CyconFKS-87}.

Clearly, if $u\in D(H)$ then $\lambda$ is an eigenvalue. But much less
restrictive growth conditions suffice to construct a Weyl sequence
from $u$ by a cut-off procedure. One of the main objectives of the
present paper is to provide a proof along these lines for a great
variety of operators. In our framework, the principal part $H_0$ of
$H$ is the selfadjoint operator associated with a strongly local
regular Dirichlet form $\ee$ and $H=H_0+\mu$ with a measure
perturbation. Of course, this includes Schr\"odinger operators on
manifolds and open subsets of Euclidean space, but much more singular
coefficients are included. In our general Sch'nol's theorem potentials
in $L^1_{\text{loc}}$ with form small negative and arbitrary positive
part are included, thereby generalizing results that require some Kato
class condition. The appropriate ``decay assumption'' on $u$ that is
necessary can roughly be called subexponential growth and is phrased
in terms of conditions like
$$
\frac{\| u\chi_{B(x_0,r_n+\delta)}\|}{\| u\chi_{B(x_0,r_n)}\|}\to 1\mbox{
  for some }r_n\to\infty
$$
and some fixed $\delta >0$.  Here, $\chi_M$ is the characteristic function of $M$ and  $B(p,s)$ denotes the closed ball in the intrinsic metric around $p$ with radius $s$. A precise definition of the intrinsic metric is given below.  
%The right scale is given by the intrinsic
%metric, with which these balls are defined. 
For uniformly bounded and
strictly elliptic divergence form operators, one recovers the usual
Euclidean balls.

It is interesting to note that we use a form analog of Weyl sequences
that enables us to treat partial differential operators with singular
coefficients. Of course, the usual calculations of $H(\eta u)$ for a
smooth cut-off function $\eta$ fail in the present general context. That is already
true for operators in divergence form with nondifferentiable coefficients and to our 
knowledge, there is no Sch'nol's Theorem in that context available in the literature.  
They have to be replaced by calculations with the corresponding forms.
The crucial object in that respect is the \emph{energy measure} of a
strongly local Dirichlet forms that supplies one with a calculus
reminiscent of gradients. All this together leads to our version of
Sch'nol's theorem, Theorem \ref{SchnolConcrete} below which is one of
the main results of the present paper. Apart from its generality it is
also pretty simple conceptually.

Another aim of the present paper is to advertise Dirichlet form
techniques for \emph{quantum} or \emph{metric graphs}. As a space
these consist of a countable family of edges (intervals) that are
glued together in the sense that the Laplacian on the direct sum of
intervals is equipped with certain boundary conditions for those edges
that meet at a vertex. For certain types of boundary conditions one
can apply the Dirichlet form framework. In this way we get a similar
understanding (and a partial generalization) of results by P. Kuchment
\cite{Kuchment-05} on Sch'nol's theorem for quantum graphs. Needless
to say that on the other hand quantum graphs provide a wealth of
examples of strongly local Dirichlet forms. While Sch'nol's theorem
had already been known for quantum graphs, the way to interpret them
as Dirichlet forms opens a powerful arsenal of analytic and
probabilistic techniques. Quite a number of results in operator and
perturbation theory have been established in the Dirichlet form
setting and can readily be applied to quantum graphs.

This can be illustrated by the ``converse'' of Sch'nol's theorem.
Proving results on ``expansion in generalized eigenfunctions'' one
gets the fact that for spectrally almost $\lambda\in\sigma(H)$ there
exists a solution that doesn't increase too seriously. In the context
of Dirichlet forms that has been established in
\cite{BoutetdeMonvelS-03b}; see also the references in there and the
discussion in \cite{CyconFKS-87}. Together with what we said above,
the results from \cite{BoutetdeMonvelS-03b} can directly be applied to
certain quantum graphs which yields a partial converse of Kuchment's
results in \cite{Kuchment-05} that seems to be new.

At least in terms of existing proofs the ``Sch'nol part'' of the
characterization of the spectrum in terms of eigenfunctions appears to
be the easier one. That is reflected in the fact that we needed more
restrictive conditions in \cite{BoutetdeMonvelS-03b} to establish an
eigenfunction expansion than what we need in the present paper. That
refers to conditions on the underlying operator as well as to
conditions on the measure perturbation, where a Kato type condition is
needed in \cite{BoutetdeMonvelS-03b}. The conclusion from the latter
paper is that for spectrally almost $\lambda\in\sigma(H)$ there is a
``subexponentially bounded'' eigensolution. To see that this is
compatible with the growth condition referred to above is the third
main result of the present paper. We should, moreover, mention our
version of the Caccioppoli inequality, Theorem \ref{Caccioppoli}
below. For the unperturbed operator $H_0$ such an inequality can be
found in \cite{BiroliM-95}. Our version here, including measure
perturbations, appears to be new and might be of interest in its own
right.

\section{Assumptions and basic properties}
\subsection*{Dirichlet forms}

Throughout we will work with a locally compact, separable metric space
$X$ endowed with a positive Radon measure $m$ with $ \mbox{supp} m=X$.
Our exposition here goes pretty much along the same lines as those in
\cite{BiroliM-95,Sturm-94b}. We refer to \cite{Fukushima-80} as the
classical standard reference as well as
\cite{BouleauH-91,FukushimaOT-94,MaR-92,Davies-90b} for literature on
Dirichlet forms.  The central object of our studies is a regular
Dirichlet form $\mathcal{E}$ with domain $\mathcal{D}$ in $L^2(X)$ and
the selfadjoint operator $H_0$ associated with $\mathcal{E}$.  This
means that $\mathcal{D} \subset L^2(X,m)$ is a dense subspace,
$\mathcal{E}\colon\mathcal{D} \times \mathcal{D} \rightarrow \KK$ is
sesquilinear and $\mathcal{D}$ is closed with respect to the energy
norm $\|\,\cdot\,\|_\mathcal{E}$, given by
\[
\|u\|_\mathcal{E}^2=\mathcal{E}(u,u) +\| u\|_{L^2(X,m)}^2,
\]
in which case one speaks of a \textit{closed form} in $L^2(X,m)$. In the sequel we will write
$$\mathcal{E}(u):= \mathcal{E} (u,u). $$
Let us emphasize that in contrast to most of the work done on
Dirichlet forms, we explicitly include the case of complex scalars;
$\KK$ denotes either $\RR$ or $\CC$.  The unique operator $H_0$
associated with $\mathcal{E}$ is then characterized by
\[
D(H_0)\subset \mathcal{D} \ \mbox{and } \mathcal{E}(f,v)=(H_0f\mid
v)\quad (f\in D(H_0), v\in \mathcal{D}).
\]
Such a closed form is said to be a \textit{Dirichlet form} if
$\mathcal{D}$ is stable under certain pointwise operations; more
precisely, $T:\KK\to\KK$ is called a \emph{normal contraction} if
$T(0)=0$ and $|T(\xi)-T(\zeta)|\le |\xi -\zeta|$ for any
$\xi,\zeta\in\KK$ and we require that for any $u\in \mathcal{D}$ also
\[
T\circ u\in \mathcal{D}\mbox{ and }\mathcal{E}(T\circ u)\le
\mathcal{E}(u).
\]
Here we used the original condition from \cite{BeurlingD-58} that
applies in the real and the complex case at the same time. Today,
particularly in the real case, it is mostly expressed in an equivalent
but formally weaker statement involving $u\vee 0$ and $u\wedge 1$, see  
\cite{Fukushima-80}, Thm. 1.4.1 and  \cite{MaR-92}, Section I.4.

A Dirichlet form is called \textit{regular} if $\mathcal{D} \cap
C_c(X)$ is dense both in $(\mathcal{D},\| \cdot \|_\mathcal{E})$ and
$(C_c(X),\| \cdot \|_{\infty })$, where $C_c(X)$ denotes the space of
continuous functions with compact support.

\subsection*{Strong locality and the energy measure}
$\mathcal{E}$ is called \textit{strongly local} if
$$\mathcal{E}(u,v)=0$$ 
whenever $u$ is constant a.s. on the support of $v$.

The typical example one should keep in mind is the Laplacian
$$H_0=-\Delta \mbox{ on }L^2(\Omega ),\quad \Omega \subset \mathbb{R}^d\mbox{ open, }$$
in which case
$$
\mathcal{D}=W^{1,2}_0(\Omega )\mbox{ and
}\mathcal{E}(u,v)=\int_{\Omega }(\nabla u| \nabla v)dx.$$ Now we turn
to an important notion generalizing the measure $(\nabla u |\nabla
v)dx$ appearing above.

In fact, every strongly local, regular Dirichlet form $\mathcal{E}$
can be represented in the form
\[
\mathcal{E}(u,v) = \int_X d\Gamma (u,v)
\]
where $\Gamma $ is a nonnegative sesquilinear mapping from
$\mathcal{D}\times\mathcal{D}$ to the set of $\KK$-valued Radon
measures on $X$. It is determined by
\[
\int_X \phi\, d\Gamma(u,u) = \mathcal{E}(u,\phi u) -\frac12
\mathcal{E}(u^2,\phi )
\]
and called \textit{energy measure}; see also \cite{BouleauH-91}.  The
energy measure satisfies the Leibniz rule,
\[
d\Gamma(u\cdot v,w)=ud\Gamma (v,w)+vd\Gamma (u,w),
\]
as well as the chain rule
\[
d\Gamma (\eta (u),w)=\eta'(u)d\Gamma (u,w).
\]
One can even insert functions from $\mathcal{D}_{\text{loc}}$ into
$d\Gamma$, where
$$\mathcal{D}_{\text{loc}} := \lbrace u\in
L^2_{\text{loc}}\mbox{ such that }\phi u\in \mathcal{D}\mbox{ for all
}\phi \in \mathcal{D}\cap C_c(X)\} ,$$ as is readily seen from the
following important property of the energy measure, \textbf{strong
  locality}:

Let $U$ be an open set in $X$ on which the function $\eta\in
\mathcal{D}_{\text{loc}}$ is constant, then
\begin{equation}\label{local}
  \chi_U d\Gamma (\eta,u) = 0,
\end{equation}
for any $u\in \mathcal{D}$. This, in turn, is a consequence of the
strong locality of $\mathcal{E}$ and in fact equivalent to the
validity of the Leibniz rule.

We write $d\Gamma(u):=d\Gamma(u,u)$ and note that the energy measure
satisfies the \textbf{Cauchy-Schwarz inequality}:
\begin{eqnarray*}
  \int_X|fg|d|\Gamma(u,v)| & \le &    
  \left(\int_X|f|^2d\Gamma(u)\right)^2\left(\int_X|g|^2d\Gamma(v)\right)^2\\
  & \le & \frac12 \int_X|f|^2d\Gamma(u)+ \frac12\int_X|g|^2d\Gamma(v) .
\end{eqnarray*}

\subsection*{The intrinsic metric}

Using the energy measure one can define the \textit{intrinsic metric}
$\rho$ by
\[
\rho (x,y)=\sup \lbrace |u(x)-u(y)|\ | u\in
\mathcal{D}_{\text{loc}}\cap C(X) \mbox{ and } d\Gamma (u)\leq
dm\rbrace
\]
where the latter condition signifies that $\Gamma (u)$ is absolutely
continuous with respect to $m$ and the Radon-Nikodym derivative is
bounded by $1$ on $X$. Note that, in general, $\rho$ need not be a
metric. (See the Appendix for a discussion of the finiteness of the $\sup$.) However, here we will mostly rely on the following

\begin{assumption}\label{A0}
  The intrinsic metric $\rho$ induces the original topology on $X$.
\end{assumption}
We denote the intrinsic balls by
$$B(x,r):=\{ y\in X| \rho(x,y)\le r\} .$$
An important consequence of the latter assumption is that the distance
function $\rho_x(\cdot):=\rho(x,\cdot)$ itself is a function in
$\mathcal{D}_{\text{loc}}$ with $d\Gamma(\rho_x)\le dm$, see
\cite{Sturm-94b}. This easily extends to the fact that for every
closed $E\subset X$ the function $\rho_E(x):= \inf\{ \rho(x,y)|y\in
E\}$ enjoys the same properties (see the Appendix). This has a very
important consequence. Whenever $\zeta : \RR\longrightarrow \RR$ is
continuously differentiable, and $\eta :=\zeta \circ \rho_E$, then
$\eta$ belongs to $\mathcal{D}_{\text{loc}}$ and satisfies
\begin{equation}\label{ac}
  d \Gamma (\eta) = (\zeta'\circ \rho_E)^2 d \Gamma (\rho_E)\leq (\zeta'\circ \rho_E)^2  dm.
\end{equation}

\subsection*{Measure perturbations}

We will be dealing with Schr\"odinger type operators, i.e.,
perturbations $H=H_0+V$ for suitable potentials $V$. In fact, we can
even include measures as potentials. Here, we follow the approach from   \cite{Stollmann-92,StollmannV-96}. Measure perturbations have been regarded 
by a number of authors in different contexts, see e.g.  
\cite{AlbeverioM-91,Hansen-99,Sturm-92} and the references there.
To set up the framework, we first recall that every regular Dirichlet
form $\ee$ defines a set function, the \emph{capacity}, in the following
way:
$$
\capa(U):= \inf\{ \ee(\phi)+\|\phi\|^2| \phi\in \mathcal{D}\cap
C_c(X), \phi\ge \chi_U\}
$$
for open $U$ and
$$
\capa(B):= \inf\{ \capa(U)| B\subset U, U\mbox{ open}\} .
$$
It is clear that the capacity of a set $B$ is bounded below by its
measure $m(B)$. In most cases of interest, the capacity is larger and
allows a finer distinction of sets. E.g., for the classical Dirichlet
form in one dimension, even a single point has positive capacity. We
say that a property holds quasi-everywhere, q.e. for short, if it
holds outside a set of capacity zero. We call a function $g$
\emph{quasi-continuous} if, for every $\varepsilon >0$ there is an
open set $U\subset X$ of capacity at most $\varepsilon$ such that $g$
is continuous on the complement $X\setminus U$. Every element $u\in
\mathcal{D}$ admits a quasi-continuous representative $\tilde{u}$.
Most of the times we will be sloppy in our notation and just identify
$u$ with a quasi-continuous representative.

We denote by $\mathcal{M}_0$ the set of nonnegative measures
$\mu:\mathcal{B}\to [0,\infty]$ that do not charge sets of capacity
$0$, i.e., those measures with $\mu(B)=0$ for every set $B$ with
$\capa(B)=0$. Here $\mathcal{B}$ denotes the Borel subsets of $X$ and
we stress the fact that we do not assume our measures to be locally
finite. Besides examples of the form $Vdm$, where $V$ is nonnegative
and measurable we should also mention the measure $\infty_B$, for a
given $B\subset X$, defined by $\infty_B(M)=\infty\cdot\capa(B\cap M)$
with the usual convention $\infty\cdot 0=0$.  For such a measure
$\mu_+\in \mathcal{M}_0$,
$$D(\ee +\mu_+):= \{ u\in \mathcal{D}|  \tilde{u}\in L^2(X,\mu_+)\},$$
$$
(\ee +\mu_+)(u,v):=
\ee(u,v)+\int_X\tilde{u}\overline{\tilde{v}}d\mu_+$$ defines a closed
form (that is not necessarily densely defined). We will use the
notation $\mu_+(u,v)$ for the integral in the above formula. It is
well defined since quasi-continuous versions of the same element in
$\mathcal{D}$ agree q.e. and so give the same integrals as the measure
does not charge sets of capacity zero. The selfadjoint operator on the
closure (in $L^2(X,dm)$ ) of $D(\ee +\mu_+)$ associated with the form
$\ee +\mu_+$ is denoted by $H_0+\mu_+$. A little more restriction is
needed for negative perturbations. We call $\mu_-$ \emph{admissible},
if
$$
D(\ee -\mu_-):= \{ u\in \mathcal{D}| \tilde{u}\in L^2(X,\mu_-)\},$$
$$
(\ee -\mu_-)(u,v):= \ee(u,v)-\mu_-(u,v)$$ defines a semibounded closed
form. Note that this implies that $\mu_-$ is a Radon measure in the
sense that it is finite on relatively compact sets. For an admissible
$\mu_-$ and $\mu_+\in \mathcal{M}_0$ we can define $\ee +\mu_+-\mu_-$
and the associated operator $H_0+\mu_+-\mu_-$ in the obvious way. To
get better properties of these operators we sometimes have to rely
upon more restrictive assumptions concerning the negative part $\mu_-$
of our measure perturbation. We write $\mathcal{M}_1$ for those
measures $\mu$ that are $\ee$-bounded with bound less than one; i.e.
measures for which there is a $\kappa<1$ and a $c_\kappa$ such that
$$
\mu(u)\le \kappa\ee(u) + c_\kappa\| u\|^2 .
$$
By the KLMN theorem (see \cite{ReedS-75}, p. 167) these measures are admissible. 
An important class
with very nice properties of the associated operators is the
\textbf{Kato class} and the extended Kato class. In the present framework it can be defined in the
following way: For $\mu\in \mathcal{M}_0$ and $\alpha >0$ we set
$$
\Phi(\mu,\alpha): C_c(X)_+\to [0,\infty],$$
$$
\Phi(\mu,\alpha)\varphi := \int_X\left(
  (H_0+\alpha)^{-1}\varphi\right){\tilde{ }}\quad d\mu .
$$
The extended Kato class is defined as
$$
\hat{\mathcal{S}}_K:= \{\mu\in \mathcal{M}_0|\exists \alpha>0:
\Phi(\mu,\alpha)\in L^1(X,m)'\}
$$
and, for $\mu\in \hat{\mathcal{S}}_K$ and $\alpha >0$,
$$
c_\alpha(\mu):=\| \Phi(\mu,\alpha)\|_{L^\infty(X,m)}(= \|
\Phi(\mu,\alpha)\|_{L^1(X,m)'}), c_{Kato}(\mu):=\inf_{\alpha>0} c_{\alpha}(\mu)
.
$$
The \textit{Kato class} is originally defined via the fundamental solution of the
Laplace equation in the classical case. In our setting it consists of those measures $\mu$ with $c_{Kato}(\mu)=0$.

\subsection*{Generalized eigenfunctions} As usual an element $u\in \mathcal{D}_{\mbox{loc}}$ is called a generalized eigenfunction or weak solution to the eigenvalue $\lambda$ if
$$\mathcal{E}(u,v) + \mu (u \overline{v}) = \lambda (u,v)$$
for all $v\in \mathcal{D}$ with compact support.

\section{A Weyl type criterion}\label{Weyl}
We include the following criterion for completeness. It is taken from
\cite{Stollmann-01}, Lemma 1.4.4.

\begin{prop}\label{PropWeyl} Let $h$ be a closed, semibounded form and
  $H$ the associated selfadjoint operator. Then the following
  assertions are equivalent:

  \begin{itemize}

  \item[(i)] $\lambda \in \sigma(H)$.

  \item[(ii)] There exists a sequence $(u_n)$ in $\mathcal{D}(h)$ with
    $\|u_n\|\to1$ and
$$ \sup_{v\in \mathcal{D}(h),\|v\|_h \leq 1} |(h -\lambda) [u_n,v ]| \to 0,$$
for $n\to \infty$.

\end{itemize}

\end{prop}

\begin{proof}
  $(i)\Longrightarrow (ii)$: Choose a Weyl type sequence $(u_n)$ if
  $\lambda\in \sigma_{ess} (H)$ and $u_n = u$ if there is a normalized
  eigenvector $u\in D(H)$.

  $(ii)\Longrightarrow (i)$: This is proven by contradiction. Assume
  $\lambda\in\rho (H)$. Then,
$$\sup_{n\in \N} \|(H - \lambda)^{-1} u_n\|_h =:C <\infty.$$
Therefore,
$$\|u_n\|^2 = | (h - \lambda) [u_n, (H-\lambda)^{-1} u_n]| \leq  C \sup_{v\in \mathcal{D}(h), \|v\|_h \leq 1} |(h-\lambda)[u_n,v]|$$
and the latter term tends to zero for $n\to \infty$ by assumption.
\end{proof}

We will produce a suitable sequence $(u_n)$ as above by a suitable
cutoff of generalized eigenfunctions. Note that to this end it is very
convenient that we do not have to construct elements of the operator
domain $D(H)$, a task that seems almost hopeless in the generality of
forms we are aiming at. In fact, already  for divergence form operators with singular coefficients there is no explicit description of the operator domain and the above criterion is of use in this important special case.

\section{A Caccioppoli type inequality}
In this section we prove a bound on the energy measure of a
generalized eigenfunction on a set in terms of bounds on the
eigenfunction on certain neighborhood of the set.

We need the following notation: For $E\in X$ and $b>0$ we define the
$b$-neighborhood of $E$ as
$$ B_b (E) :=\{y\in X : \rho(y,E)\leq b\}.$$

\begin{theorem}\label{Caccioppoli}
  Let $\mathcal{E}$ be a strongly local regular Dirichlet form
  satisfying Assumption {\rm \ref{A0}}.  Let $\mu_+\in \mathcal{M}_0$ and
  $\mu_- \in \mathcal{M}_1$ be given. Let $\lambda_0\in \RR$ be given.
  Then, there exists a $C = C (\lambda_0,\mu_-)$ such that for any
  generalized eigenfunctions $u$ to an eigenvalue $\lambda\leq
  \lambda_0$ of $H_0 + \mu$ the inequality
$$ \int_E d\Gamma (u) \leq \frac{C}{b^2} \int_{B_b (E)} |u|^2
dm$$ holds for any closed $E\subset X$ and any $b>0$.
\end{theorem}

\begin{remark} If it were not for the ``potential'' $\mu$, we could
  replace the neigbourhood by a collar around the boundary of $E$ (as
  will be clear from the proof). The Caccioppoli inequality replaces the familiar
commutator estimates that are used for Schr\"odinger operators.
\end{remark}

We give a proof of the theorem at the end of this section after two
auxiliary propositions.

\begin{prop}\label{aux1} 
  Let $H_0 +\mu$ be given as in the theorem and $u$ a generalized
  eigenfunction to the eigenvalue $\lambda$. Let $\eta\in \mathcal{D}$, $\eta$ real valued, 
  be arbitrary. Then,
$$\int \eta^2 d\Gamma (u) = (\lambda - \mu)( |\eta u|^2) - 2 \int \eta\,
u\, d\Gamma (\eta,u).$$
\end{prop}
\begin{proof}
  A direct calculation invoking Leibniz rule and the chain rule gives
  \begin{eqnarray*}
    \int \eta^2 d\Gamma (u) &=& \int d\Gamma (u,\eta^2 u) - \int u d\Gamma
    (u,\eta^2)\\ &=& \int d \Gamma (u,\eta^2 u) - 2 \int u \eta d \Gamma
    (u,\eta)\\ &=& \mathcal{E} (u,\eta^2 u) - 2 \int u \eta d\Gamma
    (u,\eta)\\ &=& (h -\lambda) (u,\eta^2 u) + (\lambda - \mu) (|\eta u
    |^2) - 2 \int u \,\eta\, d\Gamma (u,\eta).
  \end{eqnarray*}
  As $u$ is a generalized eigenfunction, the statement follows.
\end{proof}

\begin{prop}\label{aux2} Let $u,\eta \in \mathcal{D}$, $\eta$ real valued, be given. Then,
$$\mathcal{E} (\eta u) = \int \eta^2 d \Gamma (u) + \int |u|^2 d\Gamma
(\eta) + 2 \int \eta u d\Gamma (u,\eta).$$
\end{prop}
\begin{proof} This is a direct calculation.
\end{proof}

We can now give the

\begin{proof}[Proof   of Theorem \ref{Caccioppoli}]
  Let $\omega =\rho_E $ and $\zeta : [0,\infty) \longrightarrow [0,1]$
  be continuously differentiable with $\zeta (0) =1$, $\zeta \equiv 0$
  on $[b,\infty]$ and $|\zeta' (t)| \leq \frac{2}{b}$ for all $t\in
  [0,\infty)$. Set $\eta :=\zeta \circ \omega$. Of course,
$$\int_E d\Gamma (u) \leq \int \eta^2 d\Gamma (u).$$ The main idea is
now to use the previous two propositions to estimate $\int \eta^2
d\Gamma (u)$ by terms of the form $\int |u|^2 d\Gamma (\eta)$ and then
to appeal to \eqref{ac}.

Here are the details: By assumption on $\mu_-$, there exists $q<1$ and
$C_q\geq 0$ with
$$\int \varphi^2 d\mu_- \leq q \mathcal{E} (\varphi) + C_q \|\varphi\|^2$$
for all $\varphi \in \mathcal{D}$.  As $\lambda \leq \lambda_0$, this
yields
\begin{eqnarray*}
  (\lambda - \mu)(|\eta u|^2 ) &\leq & \lambda \int \eta^2 |u|^2 dm + \int
  \eta^2 |u|^2 d\mu_-\\ &\leq & \lambda_0 \|\eta u\|^2 + q
  \mathcal{E}(\eta u) + C_q \|\eta u\|^2\\ &\leq & q \mathcal{E} (\eta
  u) + (\lambda_0 + C_q) \| \eta u\|^2.
\end{eqnarray*}
Combining this with Proposition \ref{aux1} we obtain
$$\int \eta^2 d\Gamma (u) \leq q \mathcal{E} (\eta u) + (\lambda_0 +
C_q) \|\eta u\|^2 - 2 \int \eta u d\Gamma (u,\eta).$$ Invoking
Proposition \ref{aux2}, we obtain
$$\int \eta^2 d\Gamma (u) \leq  q\int \eta^2 d\Gamma (u) + q \int |u|^2 d\Gamma (\eta) + (\lambda_0 +
C_q) \|\eta u\|^2 + 2(q-1) \int \eta u d\Gamma (u,\eta).$$ Application
of Cauchy Schwarz inequality to the last term yields
$$ \int \eta^2 d\Gamma (u) \leq  q\int \eta^2 d\Gamma (u) + q \int |u|^2 d\Gamma (\eta) + (\lambda_0 +
C_q) \|\eta u\|^2$$
  $$+ \frac{1-q}{S^2} \int |u|^2 d\Gamma (\eta) + S^2 (1-q)\int \eta^2 d\Gamma (u)$$
  for any $S>0$. Hence
$$( 1 - q - S^2 (1-q))\int \eta^2 d\Gamma (u) \leq (q + \frac{(1-q)}{S^2}) \int |u|^2 d\Gamma (\eta) + (\lambda_0 + C_q)\|\eta u\|^2$$
for any $S>0$. As $q<1$ and $S>0$ is arbitrary, the statement follows
with the help of \eqref{ac}. This finishes the proof.
\end{proof}

\section{A Sch'nol type result}

In this section, we first prove an abstract Sch'nol type result. We
need the following notation.  For $E\in X$ and $b>0$ we define the
$b$-collar of $E$ as
$$ A_b (E) :=\{y\in X : \rho(y,E)\leq b\:\;\mbox{and}\;\:   \rho(y,E^c)\leq b\}.$$

\begin{prop}\label{SchnolAbstract} Let $\mathcal{E}$ be a strongly local regular Dirichlet form
  satisfying Assumption {\rm \ref{A0}}.  Let $\mu_+\in \mathcal{M}_0$ and
  $\mu_- \in \mathcal{M}_1$ be given. Let $\lambda\in \RR$ with
  generalized eigenfunction $u$ be given. If there exists $b>0$ and a
  sequence $(E_n)$ of closed subsets of $X$ with
$$\frac{\| u \chi_{A_{3b}(E_n)} \|}{\| u \chi_{E_n}\|}\longrightarrow
0,n\longrightarrow 0,$$ then $\lambda$ belongs to $\sigma (H)$.
\end{prop}

\begin{proof} 
  Let $\zeta : [0,\infty)\longrightarrow [0,1]$ be continously
  differentiable, with $\zeta(0)=1$, $\zeta \equiv 0 $ on $[b,\infty)$ and
  $|\zeta'|\leq 2/b$. Let $\omega_n :=\rho_{E_n}$ and $\eta_n :=\zeta
  \circ \omega_n$. Let $u_n :=\eta_n^2 u / \|\eta_n^2 u\|$. We show
  that $(u_n)$ satisfies the assumption of Proposition \ref{PropWeyl}: Let
  $v\in\mathcal{E}$ be arbitrary.  A direct calculation involving
  Leibniz rule gives
$$ \int d\Gamma (\eta u, v) = \int d\Gamma (u,\eta v) + \int u d\Gamma (\eta,v) - \int \overline{v} d\Gamma (u, \eta)$$
for all $\eta \in \mathcal{D}$, which are real valued.  This yields
\begin{eqnarray*} (h -  \lambda)[u_n, v] &=& \frac{1}{\|\eta_n^2 u\|} \left(\int d\Gamma (\eta_n^2 u,v) + (\mu - \lambda) (\eta_n^2 u \overline{v})  \right) \\
  &=&  \frac{1}{\|\eta_n^2 u\|} \left( \int d\Gamma (u,\eta_n^2 v)   + \int u d\Gamma (\eta_n^2, v) \right. \\
  & & - \left. \int \overline{v} d\Gamma (u,\eta_n^2) +   (\mu - \lambda) (\eta_n^2 u \overline{v})  \right)\\
  &=&  \frac{1}{\|\eta_n^2 u\|}\left( \int u d\Gamma (\eta_n^2, v) - \int \overline{v} d\Gamma (u,\eta_n^2) \right)\\
  &=& \frac{2}{\|\eta_n^2 u\| }\left( \int u \eta_n d\Gamma (\eta_n,
    v) - \int \overline{v} \eta_n d\Gamma (u,\eta_n) \right)
\end{eqnarray*}
where we used in the previous to the last step that $u$ is a
generalized eigenfunction.  Cauchy-Schwarz now gives

\begin{eqnarray*}|(h - \lambda) [u_n,v]| &\leq& \frac{2}{\|\eta_n^2
    u\|} \left( \left (\int |u|^2 d\Gamma(\eta_n)\right)^{1/2} \left(
      \int \eta_n^2 d \Gamma
      (v)\right)^{1/2} \right. \\
  &+&\left. \left| \int \eta_n \overline{v} d\Gamma (u,\eta_n)\right|\right) .
\end{eqnarray*}

We will estimate the three terms on the right hand side.

As $\eta_n$ is constant outside of $A_{b}(E_n)$ we obtain from
locality and \eqref{ac}
$$\int |u|^2 d\Gamma (\eta_n)= \int_{ A_{2 b}(E_n) } |u|^2 d\Gamma
(\eta_n) \leq \frac{4}{b^2} \|\chi_{ A_{2 b}(E_n) } u\|^2.$$

As for the second term, due to $0\leq \eta_n\leq 1$, we easily find
$$ \int \eta_n^2 d\Gamma (v) \leq \int d\Gamma (v) = \mathcal{E}(v) = const.$$

We now come to the last term. As $\eta_n$ is constant outside of
$A_{b}(E_n)$, locality again gives

$$ | \int \eta_n v d\Gamma (u,\eta_n) | =  | \int_{ A_{2 b} (E_n)} \eta_n v  d\Gamma (u,\eta_n) |.$$

By Cauchy Schwarz this can be estimated by
$$ \left( \int_{ A_{2 b} (E_n)} \eta_n^2 d\Gamma (u)\right)^{1/2}
\left(\int_{ A_{2 b} (E_n)} v^2 d\Gamma (\eta_n)\right)^{1/2}.$$ By
\eqref{ac} we can estimate $ \int_{ A_{2 b} (E_n)} v^2 d\Gamma
(\eta_n)$ by $4/b^2 \|v\|^2$. By $0\leq \eta_n\leq 1$ and Theorem
\ref{Caccioppoli}, we can estimate
$$ \left( \int_{ A_{2 b} (E_n)} \eta_n^2 d\Gamma (u)\right)^{1/2}
\leq \left( \int_{ A_{2 b} (E_n)} d\Gamma (u)\right)^{1/2} \leq
\left(\frac{C}{b^2} \int_{ A_{3b} (E_n)} |u|^2 dm \right)^{1/2} .$$

Putting these estimates together shows that there exists $c>0$ with
$$ |(h - \lambda) [u_n,v]|  \leq  c \frac{\|u \chi_{A_{3b} (E_n)}\| } {\|\chi_{E_n} u\|}$$
for all $n\in \N$. As the right hand side tends to zero by our
assumption, so does the left hand side and Proposition \ref{PropWeyl}
gives the desired result.
\end{proof}

We will now specialize our considerations to subexponentially bounded
eigenfunctions. We start with a piece of notation and two auxiliary lemmas.

\medskip

A function $J: [0,\infty) \longrightarrow [0,\infty)$ is said to
be subexponentially bounded if for any $\alpha>0$ there exists a
$C_\alpha\geq 0$ with $J(r) \leq C_\alpha \exp(\alpha r)$ for all
$r>0$.
A function $f$ on a pseudo metric space   $(X,\rho)$ with measures $m$ is said to be subexponentially bounded if for some $x_0\in X$ and  $\omega (x) = \rho (x_0,x)$ the function $e^{-\alpha \omega} u$ belongs to  $ L^2
  (X,m)$ for any $\alpha >0$.

\begin{lemma} \label{aux3} Let $J: [0,\infty) \longrightarrow
  [0,\infty)$ be subexponentially bounded.  Let $b>0$ be
  arbitrary. Then, there exist for any $\delta>0$ arbitrary large
  $r>0$ with $J(r + b) \leq e^\delta J(r)$.
\end{lemma}

\begin{proof} Assume not. Then, there exists an $R_0\geq 0$ with
$$ J(r + b) > e^\delta J(r)$$
for all $r\geq R_0$. Induction then shows
$$ J(R_0 + n b) > e^{n\delta} J (R_0)$$
for any $n\in \N$. This gives a contradiction to the bounds on $J$  for $\alpha >0$
with $\alpha b <\delta$ and large $n$.
\end{proof}

\begin{lemma} 
  Let $(X,\rho)$ be a (pseudo)metric space, $m$ a measure on $X$, $x_0\in X$
  arbitrary and $\omega (x) = \rho (x_0,x)$, $B_r := B_r (x_0)$. Let
  $u : X\longrightarrow \CC$ be subexponentially bounded. Define $$J : [0,\infty)\longrightarrow
  [0,\infty), J(r) :=\int_{B_r} |u|^2 dm.$$ Then, $J$ is subexponentially bounded.
\end{lemma}

\begin{proof} For all $\alpha >0$, we find
  \begin{eqnarray*}
    J(r) = \int_{B_r} |u|^2 dm & = & \int_{B_r} |e^{\alpha \omega}
    e^{-\alpha \omega} u|^2 dm\\ & = & \int_{B_r} e^{2\alpha \omega}
    |e^{-\alpha \omega} u|^2 dm\\ &\leq & e^{2 \alpha r} \int_{B_r}
    |e^{-\alpha \omega} u|^2 dm\\ &\leq & \|e^{-\alpha \omega} u\|^2 e^{2
      \alpha r}.
  \end{eqnarray*}
  This proves the lemma.
\end{proof}

\begin{theorem}\label{SchnolConcrete} Let $\mathcal{E}$ be a strongly local regular Dirichlet form
  satisfying Assumption {\rm \ref{A0}}. $x_0\in X$
  arbitrary and $\omega (x) = \rho (x_0,x)$. Let $\mu_+\in \mathcal{M}_0$ and
  $\mu_- \in \mathcal{M}_1$ be given. Let $u$ be a generalized
  eigenfunction which is subexponentially growing, i.e. $e^{-\alpha \omega} u\in L^2 (X,m)$ for any $\alpha>0$. 
Then, $\lambda$
  belongs to $\sigma(H)$.
\end{theorem}
\begin{proof} As $u$ is subexponentially growing, the function
$$J(r) :=\int_{B_r} |u|^2 dm$$ is subexponentially bounded by the previous
lemma. By Lemma \ref{aux3}, we can then choose $b>0$ and find a
sequence $(r_n)$ with $r_n\to \infty$ and $J(r_n + 3 b) / J(r_n - 3b)
\longrightarrow 1$, $n\to \infty$. As $J$ is monotonously increasing
this easily gives
$$ \frac{ J(r_n + 3b) - J(r_n - 3b)}{J(r_n) }\longrightarrow
0,\:n\longrightarrow \infty.$$ Thus, $u$ satisfies the assumption of
Proposition \ref{SchnolAbstract} with $E_n = B_{r_n}$ and the
statement follows.

\end{proof}
We infer the following result from \cite{BoutetdeMonvelS-03b} in a
form suitable for our purposes here. We need the following additional
properties of the intrinsic geometry:

\begin{assumption}\label{A1} For each $t>0$ the  semigroup  $e^{-tH_0}$ gives a map from  $L^2(X)$ to $L^\infty (X)$
  and all intrinsic balls have finite volume with subexponential
  growth:
$$
e^{-\alpha\cdot R}m(B(x,R))\to 0\mbox{ as }R\to\infty\mbox{ for all
}x\in X, \alpha>0 .
$$
\end{assumption}
With this assumption, Corollary 3.1 from \cite{BoutetdeMonvelS-03b} gives:
\begin{theorem}\label{T46}
  Let $\mathcal{E}$ be a strongly local regular Dirichlet form
  satisfying Assumptions {\rm \ref{A0}} and {\rm \ref{A1}}. Let $\mu=\mu_+-\mu_-$
  with $\mu_+\in\C{M}_0$ and $\mu_-\in \hat{\mathcal{S}}_K$ with $
  c_{Kato}(\mu)<1$. Define $H:=H_0+\mu$. Then for spectrally a.e.
  $\lambda\in \sigma(H)$ there is a subexponentially bounded
  generalized eigenfunction $u\not= 0$ with $Hu=\lambda u$.
\end{theorem}
Thus, together with Theorem \ref{SchnolConcrete}, we get the following
characterization of the spectrum:
\begin{coro}
  Let $\mathcal{E}$ be a strongly local regular Dirichlet form
  satisfying Assumptions {\rm \ref{A0}} and {\rm \ref{A1}}. Let $\mu=\mu_+-\mu_-$
  with $\mu_+\in\C{M}_0$ and $\mu_-\in \hat{\mathcal{S}}_K$ with $
  c_{Kato}(\mu)<1$. Define $H:=H_0+\mu$. Then the spectral measures of
  $H$ are supported on
$$
\{ \lambda\in\RR | \exists \mbox{ subexponentially bounded }u\mbox{
  with }Hu=\lambda u\} .
$$
\end{coro}

\section{Application: Metric and Quantum Graphs}

We now introduce a class of examples that has attracted considerable
interest in the physics as well as in the mathematical literature. We
refer the reader to \cite{AizenmanSW-06b,HislopP,Kuchment-02,
  Kuchment-04,Kuchment-05,Exner-97,KostrykinS-99b,KostrykinS-00b,KostrykinS-06}
and the references in there. Although different levels of generality
and very different ways of notation can be found in the literature,
the basic idea is the same: a \emph{metric graph} consists of line
segments -- edges -- that are glued together at vertices. In contrast
to combinatorial graphs, these line segments are taken seriously and
in fact one is interested in the Laplacian on the union of the line
segments. To get a self adjoint operator one has to specify boundary
conditions at the vertices. More precisely, we work with the following
\begin{definition}
  A \emph{metric graph} is $\Gamma=(E,V,i,j)$ where
  \begin{itemize}
  \item $E$ (edges) is a countable family of open intervals $(0,l(e))$
    and $V$ (vertices) is a countable set.
  \item $i:E\to V$ defines the initial point of an edge and $j:\{ e\in
    E| l(e)<\infty \}\to V$ the end point for edges of finite length.
  \end{itemize}
  We let $X_e:= \{ e\}\times e$, $X=X_\Gamma=V\cup\bigcup_{e\in E}X_e$
  and $\overline{X_e}:=X_e\cup\{ i(e),j(e)\}$
\end{definition}
Note that $X_e$ is basically just the interval $(0,l(e))$, the first
component is added to force the $X_e$'s to be mutually disjoint. The
topology on $X$ will be such that the mapping $\pi_e:X_e\to (0,l(e)),
(e,t)\mapsto t$ extends to a homeomorphism again denoted by
$\pi_e:\overline{X_e}\to \overline{(0,l(e))}$ that satisfies
$\pi_e(i(e))=0$ and $\pi_e(j(e))=l(e)$ (the latter in case that
$l(e)<\infty$).  To define a metric structure on $X$ we proceed as
follows: we say that $p\in X^N$ is a \emph{good polygon} if, for every
$k\in\{ 1, ... ,N\}$ there is a unique edge $e\in E$ such that
$\{x_k,x_{k+1}\}\subset \overline{X_e}$.  Using the usual distance in
$[0,l(e)]$ we get a distance $d$ on $\overline{X_e}$ and define
$$L(p)=\sum_{k=1}^Nd(x_k,x_{k+1}) .$$
Since multiple edges are, obviously, allowed, we needed to restrict
our attention to good polygons to exclude the case that
$\{x_k,x_{k+1}\}$ are joined by edges of different length.  Provided the graph is connected and
that the degree of every vertex $v\in V$
$$
d_v:= |\{ e\in E |v\in\{ i(e),j(e)\}| <\infty ,
$$
a metric on $X$ is given by
$$
d(x,y):=\inf\{ L(p)| p\mbox{ a good polygon with } x_0=x\mbox{ and
}x_N=y\} .
$$

In fact, symmetry and triangle inequality are evident and the
separation of points follows from the finiteness. Clearly, with the
topology induced by that metric, $X$ is a locally compact, separable
metric space. Note that in our setting we do allow loops, multiple
edges and there is no on upper or lower bounds for the
length of edges. In that respect, we allow more general graphs than
those considered in the literature so far. To be able to use the
framework of regular Dirichlet forms, we restrict our attention to
certain boundary conditions, known as \emph{Kirchhoff} and
$\delta$-b.c. The operator with Kirchhoff b.c. is defined as the
operator corresponding to the form
$$
\Dee=\Dee(\ee):=W^{1,2}_0(X), \:\;\ee(u,v):=\sum_e(u_e'|v_e') ,
$$
where $u_e:=u\circ \pi_e^{-1}$ defined on $(0,l(e))$,
\begin{eqnarray*}
  W^{1,2}(X)&=& \left\{ u\in C(X)| 
    \sum_{e\in E}\| u_e\|_{W^{1,2}}^2=:\| u \|_{W^{1,2}}^2<\infty \right\} ,
  \\ W^{1,2}_0(X)&:=&W^{1,2}(X)
  \cap C_0(X) .\end{eqnarray*}
Clearly, $\ee$ is a regular Dirichlet form in $L^2(X,m)$, where $m$ is the 
measure induced
by the image of the Lebesgue measure on each $X_e$, so that $L^2(X,m)\ni 
u\mapsto 
(u_e)_{e\in E}\in \bigoplus_{e\in E}L^2((0,l(e)),dt)$ is unitary.

This form is strongly local with energy measure
$$d\Gamma(u,v)=\sum_{e\in E} u_e'(\pi_e(x))
\overline{v}_e'(\pi_e(x))dm(x) .$$ 
We denote by $H_0$ the operator associated with $\ee$
Note that every point $x\in X$ has
positive capacity by the Sobolev embedding theorem so that every
measure $\mu:\C{B}\to [0,\infty]$ belongs to $\C{M}_0$.

\begin{coro} For $\ee$ as above, 
 let $\mu_+: \C{B}\to [0,\infty]$ and
  $\mu_- \in \mathcal{M}_1$ be given. Let $u\not= 0$ be a generalized
  eigenfunction  for $H:=H_0+\mu$ that is subexponentially bounded. Then, $\lambda$
  belongs to $\sigma(H)$.
\end{coro}
\begin{remark}
As we mentioned above, $\mu_+$ may include arbitrary sums of $\delta$-measures, in particular $\delta$-measures at points of $V$ for which one gets a quantum graph with $\delta$-boundary conditions with positive coefficients
\end{remark}
For an application of Theorem \ref{T46} we have to require more restrictive
conditions, which, however, are met in many examples. 
% \begin{remark}\begin{enumerate}
%\item
%We have
%$$
%\| (H_0+\alpha)^{-1}:L^1(X)\to C_0(X)\|\to 0\mbox{  as }\alpha\to\infty .
%$$
%In fact, the Sobolev imbedding theorem gives $(H_0+\alpha)^{-\frac12}:L^2(X)\to C_0(X)$ so that $(H_0+\alpha)^{-\frac12}:L^1(X)\to L^2(X)$ by duality. In particular, 
%$$
%e^{-tH_0}:L^2(X)\to L^\infty (X)
%$$
%so that the first part of \ref{A1} is met.
%\item The Combes-Thomas estimate (see \cite{BoutetdeMonvelS-03b} for a version in the Dirichlet form setting) implies that
%$$
%\| \chi_A (H_0+\alpha)^{-1} \chi_B\| \le \exp(-\gamma_\alpha\dist(A,B) $$
%where $\gamma_\alpha\to\infty$ for $\alpha\to\infty$.
%\item 
%From (1) and (2) we get that any measure that is exponentially bounded in
%the sense that there is a $\beta>0$ with 
%$$
%\sup_{x\in X, R>0}e^{-\beta R}\mu (B(x,R))< \infty
%$$
%is in the Kato class with $c_{Kato}(\mu)=0$.
%\end{enumerate}
%\end{remark}
%The above remarks show that the following result applies in great generality:

\begin{coro}
Let $\ee$ be the Dirichlet form of a metric graph $\Gamma$ as above. Assume that {\rm \ref{A1}} is satisfied. Let $\mu =\mu_+-\mu_-$ where
$\mu_-\in \hat{\mathcal{S}}_K$ with $
  c_{Kato}(\mu)<1$. Define $H:=H_0+\mu$. Then for spectrally a.e. $\lambda\in
\sigma(H)$ there exists a subexponentially bounded $u\not= 0$  with $Hu=\lambda u $.
\end{coro}

For certain tree graphs an expansion in generalized eigenfunctions has been given in \cite{HislopP}. In a forthcoming work 
\cite{LSSV} we will prove that generalized eigenfunction expansions exist for much more general graphs than treated above.

\begin{appendix}
\section{Properties of  absolutely continuous 
elements, the distance function $\rho_E$ and 
all that}
Let $\mathcal{E}$ be a regular strongly local Dirichlet form with
associated energy measure $\Gamma$.  In this appendix, we discuss some
properties of
$$\mathcal{A}:=\{u\in \mathcal{D}_{\text{loc}} : u\;\mbox{real valued with}\;\:d\Gamma (u) \leq
dm\}.$$ We apply this to show that $\rho_E$ belongs to $\mathcal{A}$
for any closed $E\subset X$ (and in fact for any $E\subset X$) if
\ref{A0} is statisfied.  For $E$ consisting of single points this was
first shown in \cite{Sturm-94b}. For closed $E$ this seems to be
known. It is stated for example in \cite{Sturm-95b}, where a proof is
attributed to \cite{Sturm-94b}. As we did not find the proof there, we
could not resist to produce one here. Along our way we will also
reprove the case of a single point. Moreover, we will discuss
connectedness of the space $X$ in terms of the intrinsic metric.

\bigskip

We start by collecting basic properties of $\mathcal{A}$.
\begin{prop}\label{BasicPropertiesA}
  \begin{itemize}
  \item[(a)] $\mathcal{A}$ is balanced, i. e.  convex and closed under
    multiplication by $(-1)$.
  \item[(b)] $\mathcal{A}$ is closed under taking minima and maxima.
  \item [(c)] $\mathcal{A}$ is closed under adding constants.
  \item[(c)] $\mathcal{A}$ is closed under pointwise convergence of
    functions, which are uniformly bounded on compact sets.
  \end{itemize}
\end{prop}
\begin{proof}
  (a) Obviously, $\mathcal{A}$ is closed under multiplication by $-1$.
  Let $u,v\in \mathcal{A}$ and $\lambda, \mu \geq 0$ with $\mu +
  \lambda =1$ be given. Set $w = \lambda u + \mu v$. Then, for every
  $\varphi \in C_c (X)$ we have
  \begin{eqnarray*}
    \int \varphi^2  d\Gamma (w) &=& \lambda^2 \int
    \varphi^2 d\Gamma (u) + 2 \lambda \mu \int \varphi^2 d\Gamma (u,v) +
    \mu \int \varphi^2 d\Gamma (v)\\  &\leq&  \lambda^2 \int \varphi^2
    d\Gamma (u) + 2 \lambda \mu \left(\int \varphi^2 d\Gamma
      (u)   \int \varphi^2 d\Gamma(v) \right)^{1/2} + \mu \int \varphi^2 d\Gamma (v)\\ &\leq& \int
    \varphi^2 dm.
  \end{eqnarray*}
  As $\varphi$ was arbitrary the statement follows.

  (b) As $\mathcal{A}$ is closed under multiplication by $-1$, it
  suffices to consider minima.  A direct consequence of locality is
  the truncation property
$$ d\Gamma (u \wedge v, w) = \chi_{\{ u < v\}} d\Gamma (u,w) + \chi_{\{u
  \geq v\}} d\Gamma (v,w)$$ for all $u,v,w\in
\mathcal{D}_{\text{loc}}$.  If $w = u\wedge v$ we obtain
$$ d\Gamma (u \wedge v, u\wedge v) = \chi_{\{ u < v\}}d\Gamma (u,u) +
\chi_{\{u\geq v\}} d\Gamma (v,v).$$ This shows that $ \mathcal{A}$ is
closed under $\wedge$.

(c) This is obvious.

(d) Let $(u_n)$ be a sequence in $\mathcal{A}$ which converges
pointwise to $u$ and is uniformly bounded on each compact set.  We
first show that $u$ belongs to $\mathcal{D}_{\text{loc}}$.  Let $\psi
\in C_c (X)\cap \mathcal{D}$ be arbitrary.  Leibniz rule, Cauchy
Schwarz inequality and locality of $d\Gamma$ give
\begin{eqnarray*}\mathcal{E} (\psi u_n) &= &\int d\Gamma (\psi u_n)\\
  & =& \int
  \psi^2 d\Gamma (u_n) + 2 \int \psi u_n d\Gamma (u_n,\psi) + \int u_n^2
  d\Gamma (\psi)\\
  &\leq & \int
  \psi^2 d\Gamma (u_n) + \int u_n^2 d\Gamma (\psi)  + \int \psi^2 d\Gamma (u_n)  + \int u_n^2 d\Gamma (\psi). \\
  &\leq & 2 \int \psi^2 dm + 2\int u_n^2 \chi_{\supp \psi} d\Gamma (\psi).
\end{eqnarray*}

The assumptions on $(u_n)$ show that $(\mathcal{E} (\psi u_n))$
remains bounded. By semicontinuity of $\mathcal{E}$ we infer that
$\psi u$ belongs to $\mathcal{D}$. As $\psi \in \mathcal{D} \cap C_c
(X)$ is arbitrary, we obtain $u\in \mathcal{D}_{\text{loc}}$.

Let now an arbitrary $\varphi \geq 0$ continuous with compact support
be given. Choose $\psi$ in $\psi \in \mathcal{D} \cap C_c (X)$ with
$\psi \equiv 1$ on the support of $\varphi$.  This is possible as
$\mathcal{E}$ is a Dirichlet form.  Then, by Banach/Saks theorem,
boundedness of $(\mathcal{E} (\psi u_n))$ implies convergence of
convex combinations $(w_k)$ of the $(\psi u_n)$ with respect to the
energy norm. By convexity of $\mathcal{A}$, these convex combinations
have the form $w_k = \psi v_k$ with $v_k\in \mathcal{A}$.  As $\psi
u_n$ converge to $\psi u$ in $L^2$ we infer that the energy norm limit
of the $(w_k)$ is also $\psi u$.  Locality and convergence of $w_k =
\psi v_k$ to $u \psi$ with respect to the energy norm yield
$$\int \varphi d\Gamma (u) = \int \varphi d\Gamma ( \psi u) = \lim \int \varphi d\Gamma (\psi v_n )
=\lim \int \varphi d\Gamma ( v_n) \leq \int \varphi dm.$$ As
$\varphi\geq 0$ with compact support is arbitrary, the statement
follows.
\end{proof}

The previous proposition implies that $\mathcal{A}$ is also closed
under taking suitable suprema and infima. This is discussed next.

\begin{lemma} \label{closedness} Let $\mathcal{F} \subset
  \mathcal{A}\cap C(X)$ be stable under taking maxima (minima).  If $
  u:= \sup\{ v : v\in\mathcal{F}\}$ ($u:=\inf\{ v :
  v\in\mathcal{F}\}$) is continuous, then $u$ belongs to
  $\mathcal{A}$.
\end{lemma}
\begin{proof}
  By our assumptions on $X$, there exist compact $K_n \subset X$,
  $n\in \N$, with $X = \cup_{n\in \N} K_n$ and $K_n \subset
  K_{n+1}^\circ$. By (d) of the previous proposition it suffices to
  construct $u_n\in \mathcal{F}$ with $|u_n - u | \leq 1/n$ on $K_n$
  for each $n\in \N$. This will be done next: For $n\in \N$ and $x\in
  K_n$, we can find $v_{x,n}\in \mathcal{F}$ with $ u(x) -
  \frac{1}{2n} \leq v_{x,n} (x).$ By continuity of $v_{x,n}$ and $u$,
  there exists then an open neigbourhood $U_{x,n}$ of $x$ with
$$ u(y) - \frac{1}{n} \leq v_{x,n} (y)$$ for all $y\in U_{x,n}$. As
$K_n$ is compact, there exist $x_1,\ldots, x_l$ with $K_n \subset
\cup_{j=1}^l U_{x_j,n}$. As $\mathcal{F}$ is closed under taking
maxima, the function
$$ u_n:=\max\{ v_{x_j,n} : j=1,\ldots, l\}$$
belongs go $\mathcal{F}$. By construction
$$ u(x) - \frac{1}{n} \leq u_n \;\:\mbox{on}\:\; K_n.$$ As the 
inequality $u_n \leq u$ is clear, the proof is finished.
\end{proof}

We now turn to the distance function $\rho$.  By definition we have
$$\rho(x,y) :=\sup\{u(x) - u(y) : u\in\mathcal{A}\cap C(X)\}.$$ Direct
arguments show that $\rho (x,y)$ is nonnegative, symmetric and
satisfies the triangle inequality.  As $\mathcal{A} \cap C(X)$ is
closed under adding constants, for each $x\in X$, the distance
function $\rho_x (y) := \rho (x,y)$ is then given by
$$\rho_x (y) :=\sup\{ u(y): u \in \mathcal{F}_x\}$$ with $\mathcal{F}_x
= \{u \in \mathcal{A}\cap C (X) : u(x) =0\}$.
The following proposition is essentially contained in \cite{Sturm-94b}, page $191$ and page $194$. 

\begin{prop} Assume {\rm \ref{A0}}.  Let $x\in X$ be arbitrary. Then, $\{y :
  \rho_x (y)<\infty\}$ is exactly the connected component of $x$.
\end{prop}
\begin{proof} Set $ C_x:=\{y : \rho_x (y)<\infty\}$.  Of course, all
  functions which are constant on each component of $X$ belong to
  $\mathcal{A}\cap C(X)$. Thus, $\rho_x (y) = \infty$ whenever $x$ and
  $y$ belong to different components.  Thus, $C_x$ is contained in the
  connected component of $x$. We now show the reverse inclusion. To do
  so it suffices to show that $C_x$ is both open and closed.  By
  Assumption {\rm \ref{A0}} the set $C_x$ is open. Moreover, if $y$ belongs
  to $X\setminus C_x$, then by
$$ \infty = \rho(x,y) \leq \rho (y,z) + \rho (z,x)$$ 
we obtain that any $z\in X$ with $\rho(z,y) <1$ belongs to $X\setminus
C_x$ as well. By {\rm \ref{A0}} again the set of such $z$ is open, and the
complement $X\setminus C_x$ is shown to be open as well.
\end{proof}

\begin{prop} Assume {\rm \ref{A0}}.  Let $x\in X$ be arbitrary and $C_x$ be
  the connected component of $x$. Then, $ \chi_{C_x} \rho_x$ belongs to
  $\mathcal{A}\cap C(X)$.
\end{prop}
\begin{proof} It suffices to consider the case that $X$ is connected.
  By Assumption {\rm \ref{A0}} and the previous lemma, $\rho_x$ is then
  continuous. As $\mathcal{F}_x = \{u \in \mathcal{A}\cap C (X) : u(x)
  =0\}$ is closed under taking maxima and $\rho_x (y) =\sup\{ u : u
  \in \mathcal{F}_x\}$, the statement now follows from Lemma
  \ref{closedness}.
\end{proof}

We now turn to distances from arbitrary sets.  For $E\subset X$ we
define
$$ \rho_E (z):=\inf\{ \rho_x (z) : x\in E\}.$$

\begin{theorem}  
  Assume {\rm \ref{A0}}. Let $E\subset X$ be arbitrary and let $C$ be the
  union of the connected components of the points of $E$. Then, the function  $\chi_C \rho_E$  belongs to $\mathcal{A}\cap
  C(X)$.
\end{theorem}
\begin{proof} As $C$ is open and closed it suffices to consider the
  case $C=X$. By {\rm \ref{A0}} and triangle inequality, the function
  $\rho_E$ is continuous. Moreover, as discussed above $\rho_x$
  belongs to $\mathcal{A}\cap C(X)$ for any $x\in X$. The statement
  now follows from Lemma \ref{closedness}. \end{proof}

% for $F = \{x_1,\ldots,x_n\}$ we have $\rho_F = \rho_{x_1}\wedge
% \ldots \wedge \rho_{x_n}$ and (b) of the preceeding proposition
% gives
%$$d\Gamma (\rho_{F}) \leq dm.$$  
%As $X$ is separable, so is any closed subset $E\subset X$. We can
% therefore find $F_n\subset X$ finite such that $(\rho_{F_n})$
% converge pointwise to $\rho_E$ and are uniformly bounded on compact
% sets.  By (d) of the preceeding proposition, $\mathcal{A}$ is closed
% under such limits and we see that $\rho_E$ belongs to $\mathcal{A}$
% for all closed $E\subset X$.

We note a consequence of the previous theorem.
\begin{coro} Assume {\rm \ref{A0}}.  For $E\subset X$, the equality $\rho_E
  (z) =sup\{ u(z) : u\in \mathcal{F}_E\}$ holds, where $\mathcal{F}_E
  :=\{v\in \mathcal{A}\cap C(X) : v \equiv 0\,\mbox{on}\: E\}$.
\end{coro}
\begin{proof}
  Denote the supremum in the statement by $\rho_E^*$. As $\rho_x
  (z)\geq u(z)$ for any $u\in \mathcal{F}_E$ and $x\in E$, we have
  $\rho_E \geq \rho_E^*$. For the converse direction, we note that
  $\rho_E$ belongs $\mathcal{F}_E$ by the previous theorem.
\end{proof}

We finish this section by noting a strong closedness property of
$\mathcal{A}$.

\begin{prop} 
  $\mathcal{A}$ is closed under convergence in $L^2_{\text{loc}}$.
\end{prop}
\begin{proof} Let $K$ be an arbitrary compact subset of $X$.  As
  $\rho_x$ belongs to $\mathcal{A}$ for any $x\in X$, we can find
  $\psi \in C_c (X) \cap \mathcal{A}$ with $ \psi \equiv 1$ on $K$
  (take e.g $\psi :=\max\{0,\frac{1}{R} \min\{R, 2 R - \rho_x\}\}$ for
  $x\in K$ and $R$ large).  The proof follows by mimicking the
  argument in the proof of (d) Proposition \ref{BasicPropertiesA} and
  using that $d\Gamma (\psi) \leq dm$.
  % Let $(u_n)$ be a sequence in $\mathcal{A}$ which converges in
  % $L^2_{\textit{loc}}$ to $u$.  Let $K$ be an arbitrary compact
  % subset of $X$.  As $\rho_x$ belongs to $\mathcal{A}$ for any $x\in
  % X$, we can find $\psi \in C_c (X) \cap \mathcal{A}$ with $ \psi
  % \equiv 1$ on $K$ (take e.g $\psi :=\max\{0,\frac{1}{R} \min\{R, 2
  % R - \rho_x\}\}$ for $x\in K$ and $R$ large).  Let $K'$ be the
  % compact support of $\psi$. Mimicking the first part of the
  % considerations of the proof of (d) Proposition
  % \ref{BasicPropertiesA} together with $d\Gamma (\psi) \leq dm$
  % gives
%$$\mathcal{E} (\psi u_n)  \leq  2\int \psi^2 dm  + \int \chi_{K'} u_n^2 dm.$$
%This shows boundedness of $\mathcal{E} (\psi u_n)$. We can now mimic
% the second part of the proof of (d) Proposition
% \ref{BasicPropertiesA}.
\end{proof}

\end{appendix}

%\begin{thebibliography}

%\def\cprime{$'$}\def\polhk#1{\setbox0=\hbox{#1}{\ooalign{\hidewidth
%  \lower1.5ex\hbox{`}\hidewidth\crcr\unhbox0}}}

\end{document}